\documentclass[12pt]{article}
\usepackage{amsfonts,amssymb,amscd}
\usepackage{amsmath}
\numberwithin{equation}{section}
\begin{document}
\title{The Exterior Derivative - A direct approach}
\author{
Gopala Krishna Srinivasan
\footnote{e-mail address: gopal@math.iitb.ac.in}
}
\date{
 Department of  Mathematics, Indian Institute of Technology Bombay
}
\maketitle
\footnotetext{2010 Mathematics  Subject Classification 58A10}
\paragraph*{Abstract:} In this note we provide a direct approach to the most basic operator in this theory namely the exterior derivative. The 
crucial ingredient is a calculus lemma based on determinants. We maintain the view that in a first course at least this 
direct approach is preferable to the more abstract one based on characterization of the exterior derivative in terms of 
its properties. 
\section*{1. Introduction:} 
In low dimensions at least, 
differential forms made their appearance in analysis more than three centuries ago originating in the works of Euler, Lagrange, 
Clairut and others.  
 For arbitrary orders they were 
introduced by Poincar\'e and E. Cartan. The interesting historical development is available in the papers of 
Samelson \cite{samelson}.   
A significant role was played by Pfaff (thesis advisor 
of C. F. Gauss), Jacobi and many other mathematicians. 
The final form in which this topic is currently is the culmination of efforts  
a few centuries.

The exterior derivative is one of the most important ideas in the theory of 
differential manifolds leading directly to the de Rham cohomology of manifolds - 
that is to algebraic topology ! 
 It was used by Cartan in his formulation of differential geometry via moving frames \cite{spivak2} (chapter 7).

The treatment of exterior derivative of $k$-forms 
in Hicks \cite{hicks} is completely coordinate free and described as a 
$k+1$ linear form acting on the $C^{\infty}(M)$ module $\mathfrak{X}(M)$. The defining formula (stated here for 
simplicity only for $k = 2$) is:
$$
d\omega(X, Y) = X(\omega(Y)) - Y(\omega(X)) - \omega([X, Y]). 
$$
Though the ultimate goal of introducing this operator in a coordinate free manner has been reached, 
the treatment in Hicks is somewhat austere. 
This formula is also available in Chern et al.,\cite{chern} or \cite{spivak1} (p. 213) as well.  
 The proofs  in \cite{chern} and \cite{spivak1} on 
 the existence and uniqueness of the exterior differential operator (satisfying certain conditions) employs a mixture of 
local and global arguments (see also \cite{willmore}). Although the construction of the exterior derivative carried out in 
\cite{chern} (through its characterization in terms of its properties) is quite elegant we feel an alternate 
treatment which is direct would be useful 
for audience in a first course (following the books of M. Spivak \cite{spivak} or J. R. Munkres \cite{munkres}) keeping in focus certain special tensorial features.   
Specifically, the case with the exterior derivative is reminiscent of that of 
covariant derivative namely,  when changing coordinates, 
certain undesirable terms ought to cancel out.  We see this happen explicitly here in our discussion of the exterior derivative.

We show that the key ingredient needed for defining the exterior derivative 
is a basic calculus identity involving determinants. 
We find this interesting  inasmuch as 
when changing coordinates we actually witness the internal cancellations of terms involving 
the second derivatives of the transition maps. 
\section*{2. Basic Calculus Lemma and the existence of $d$:} 
In the following lemma $J$ would denote an ordered $k-$set $\{j_1, j_2, \dots j_k\}$, 
$1 \leq j_1 < j_2 < \dots < j_k \leq n$.   We shall consider pairs $(s, J)$ such that $1 \leq s \leq n$ and 
$s \neq j_1, j_2, \dots, j_k$. Let $N$ be the total number of such pairs. 
For a given pair $(s, J)$ and $j_p \in J$ denote by $(j_p, J^{\prime})$ 
the complementary pair $(j_p, J^{\prime})$ obtained by removing $j_p$ from $J$ and inserting $s$ in the right place. 
\paragraph*{Lemma:} Let $\phi_1, \phi_2,\dots, \phi_k$ be $k$ smooth functions of $z_1, z_2, \dots, z_n$. Then 
$$
\sum_{s, J}
\frac{\partial}{\partial z_s}
\frac{\partial(\phi_1, \phi_2, \dots, \phi_k)}{\partial(z_{j_1}, z_{j_2}, \dots, z_{j_k})} dz_s \wedge dz_{j_1} \wedge \dots \wedge 
dz_{j_k}= 0. 
$$
\paragraph*{Proof:} Assume $j_{q-1} < s < j_q$ so that $dz_s$ would need to move through $q-1$ transpositions
 to bring the monomial in 
standard form leading to a factor of $(-1)^{q-1}$.  Also carrying out the indicated differentiation would produce
 $k$ determinants out of each 
summand leading to $kN$ monomials in all. We need to show that the monomials
 can be paired off in such a way that the sum is ultimately zero.

Let us consider the terms coming from the complementary pairs $(s, J)$ and $(j_p, J^{\prime})$. We may assume at the outset 
that 
$s < j_p$ for in the opposite case we can interchange the roles of $(s, J)$ and $(j_p, J^{\prime})$. 
The determinant will be 
written in such a way that the second derivatives appear in the first column which necessitates a book-keeping of the number of 
column exchanges.

Computing the derivative of the determinant, the term wherein the $p$th column is differentiated is:
$$
(-1)^{p-1}\begin{vmatrix}
\frac{\partial^2\phi_1}{\partial z_s\partial z_{j_p}} & \frac{\partial \phi_1}{\partial z_{j_1}} & \dots & 
\frac{\partial \phi_1}{\partial z_{j_k}} \\
\dots & \dots & \dots &\dots\\
\frac{\partial^2\phi_1}{\partial z_s\partial z_{j_p}} & \frac{\partial \phi_1}{\partial z_{j_1}} & \dots & 
\frac{\partial \phi_1}{\partial z_{j_k}} \\
\end{vmatrix}
$$
Together with the differentials, we get the monomial:
$$
(-1)^{p+q-2}\begin{vmatrix}
\frac{\partial^2\phi_1}{\partial z_s\partial z_{j_p}} & \frac{\partial \phi_1}{\partial z_{j_1}} & \dots & 
\frac{\partial \phi_1}{\partial z_{j_k}} \\
\dots & \dots & \dots &\dots\\
\frac{\partial^2\phi_1}{\partial z_s\partial z_{j_p}} & \frac{\partial \phi_1}{\partial z_{j_1}} & \dots & 
\frac{\partial \phi_1}{\partial z_{j_k}} \\
\end{vmatrix} dz_{j_1}\wedge\dots\wedge dz_{j_{q-1}}\wedge dz_s \wedge dz_{j_q}\wedge \dots \wedge dz_{j_k} \eqno(1)
$$
Now we consider the term arising out of the complementary pair $(j_p, J^{\prime})$ and look at the relevant monomial namely 
the one in which the 
second derivatives 
$$
\frac{\partial^2\phi_i}{\partial z_{j_p}\partial z_{s}}, \quad i = 1, 2, \dots, k
$$
appear in the determinant. 
 These second derivatives appear in the $q$th column ($j_{q-1} < s < j_q$) 
and so we need $q-1$ column exchanges to bring them to  the first column 
thereby producing a $(-1)^{q-1}$ sign. We also have in addition the factor 
$$
dz_{j_p}\wedge dz_{j_1} \wedge \dots \wedge dz_{j_{q-1}} \wedge dz_s \wedge dz_{j_q} \wedge \dots \wedge dz_{p-1}   
\wedge dz_{p+1} \wedge \dots 
\wedge dz_{j_k},\quad (s < j_p).
$$
Owing to the presence of $dz_s$, the $dz_{j_p}$ has to now 
move through $p$ transpositions to get this 
in standard form. Thus we get the term (1) but with $(-1)^{p+q-1}$ instead. Thus the terms arising from complementary 
pairs cancel out. The proof is complete. 
\paragraph*{Definition (The exterior derivative $d$):} The standard notation for the set of all smooth $k$-forms on $M$ is 
$\Omega^k(M)$. We introduce the $\mathbb R$-linear map 
$$
d : \Omega^k(M) \longrightarrow \Omega^{k+1}(M). 
$$
Let $\omega$ be a differential $k$ form and on the chart $U$ let $\omega$ be given by 
$$
\omega = \sum_i a^U_{i_1i_2\dots i_k} dx_{i_1} \wedge \dots \wedge dx_{i_k}
$$
We define on each chart  
$$
d\omega = \sum_i da^U_{i_1i_2\dots i_k} \wedge dx_{i_1} \wedge \dots \wedge dx_{i_k}. \eqno(2)
$$
The notation $\displaystyle{\sum_i
}$ stands for the sum over all standard $k-$tuples $(i_1, i_2, \dots, i_k)$ with 
$i_1 < i_2 < \dots < i_k$.

The basic properties of this operator can almost be \it read off \rm  from this definition except for one hurdle. 
We do have the job of showing that $d$ is well defined namely, of 
verifying consistency on overlapping charts but we shall demonstrate that this 
is nothing but the basic calculus lemma! 
\paragraph*{Theorem:} The operator $d$ given by (2) is a well-defined element of 
$\Omega^{k+1}(M)$. 
\paragraph*{Proof:} 
Let $U$ and $V$ be two overlapping charts. 
Need to check that 
$$
\sum_i da^U_{i_1i_2\dots i_k} \wedge dx_{i_1} \wedge \dots \wedge dx_{i_k} = 
\sum_j da^V_{j_1j_2\dots j_k} \wedge dy_{j_1} \wedge \dots \wedge dy_{j_k},\quad \mbox{ on } U \cap V. 
$$
Since
$$
a^V_{j_1j_2\dots j_k} = 
\sum_i a^U_{i_1i_2\dots i_k}\frac{\partial(x_{i_1}, x_{i_2}, \dots, x_{i_k})}{\partial(y_{j_1}, y_{j_2},\dots, y_{j_k})} \eqno(3)
$$
our job is to check that 
$$
\sum_i da^U_{i_1i_2\dots i_k} \wedge dx_{i_1} \wedge \dots \wedge dx_{i_k} = 
\sum_i\sum_j d\Big(
a^U_{i_1i_2\dots i_k}\frac{\partial(x_{i_1}, x_{i_2}, \dots, x_{i_k})}{\partial(y_{j_1}, y_{j_2},\dots, y_{j_k})}
\Big)\wedge  
dy_{j_1} \wedge \dots \wedge dy_{j_k}
$$
Well, the right hand side breaks up into two sums:
$$
\sum_i\sum_j d\Big(
a^U_{i_1i_2\dots i_k}\frac{\partial(x_{i_1}, x_{i_2}, \dots, x_{i_k})}{\partial(y_{j_1}, y_{j_2},\dots, y_{j_k})}
\Big)\wedge  
dy_{j_1} \wedge \dots \wedge dy_{j_k} = I + II. 
$$
The first sum $I$ displayed below is tensorial namely, 
$$
\sum_i\sum_j 
da^U_{i_1i_2\dots i_k}\Big(\frac{\partial(x_{i_1}, x_{i_2}, \dots, x_{i_k})}{\partial(y_{j_1}, y_{j_2},\dots, y_{j_k})}
\Big) \wedge 
dy_{j_1} \wedge \dots \wedge dy_{j_k} = \sum_i 
da^U_{i_1i_2\dots i_k} \wedge dx_{i_1} \wedge \dots \wedge dx_{i_k}
$$
which is the desired result and so we must show that the second (non-tensorial) 
term $II$ is identically zero namely, 
$$
II = \sum_i\sum_j 
a^U_{i_1i_2\dots i_k}d\Big(\frac{\partial(x_{i_1}, x_{i_2}, \dots, x_{i_k})}{\partial(y_{j_1}, y_{j_2},\dots, y_{j_k})}
\Big) \wedge 
dy_{j_1} \wedge \dots \wedge dy_{j_k} = 0.
$$
Since the coefficients, apart (3), are arbitrary smooth functions, we must show that each of the pieces 
$$
\sum_j 
d\Big(\frac{\partial(x_{i_1}, x_{i_2}, \dots, x_{i_k})}{\partial(y_{j_1}, y_{j_2},\dots, y_{j_k})}
\Big) \wedge 
dy_{j_1} \wedge \dots \wedge dy_{j_k}
$$
individually vanishes. That is to say for each fixed $i_1 < i_2 < \dots < i_k$ we have 
$$
\sum_j 
\frac{\partial}{\partial y_s}\Big(\frac{\partial(x_{i_1}, x_{i_2}, \dots, x_{i_k})}{\partial(y_{j_1}, y_{j_2},\dots, y_{j_k})}
\Big) dy_s \wedge 
dy_{j_1} \wedge \dots \wedge dy_{j_k} = 0. 
$$
But this exactly the calculus lemma. The proof is complete. 

\end{document}